\documentclass[10pt]{article}

\usepackage[utf8]{inputenc}
\usepackage[T1]{fontenc}
\usepackage{amsmath}
\usepackage{amsfonts,amssymb,latexsym,multicol}
\usepackage[francais]{babel}
\usepackage{epsfig,epsf}
\newcounter{fig}
\input{cyracc.def}

\font\tencyr=wncyr10 scaled 1000
\def\cyr{\tencyr\cyracc}

\newtheorem{ques}{Question}

\newtheorem{prop}{Proposition}

\newcommand{\cad}{\text{c'est-à-dire }}
\newcommand{\expli}[1]{\quad\text{\footnotesize (#1)}}
\newcommand{\et}{\quad \text{et} \quad}

\newcommand{\implique}{\Rightarrow}

\newcommand{\eps}{\varepsilon}
\newcommand{\fhi}{\varphi}
\newcommand{\ioe}{\leqslant}
\newcommand{\soe}{\geqslant}
\newcommand{\vers}{\rightarrow}

\newcommand{\demi}{{\frac{1}{2}}}

\newcommand{\un}{{\bf 1}}

\newcommand{\Fcal}{{\mathcal F}}

\newcommand{\Rcal}{{\mathcal R}}

\newcommand{\Nat}{{\mathbb N}}

\newcommand{\Real}{{\mathbb R}}

\newcommand{\fin}{\hfill$\Box$}
\newcommand{\dem}{\noindent {\bf D\'emonstration\ }}
\newcommand{\fine}{\tag*{\mbox{$\Box$}}}

\title{Remarques élémentaires sur la fonction de Möbius}
\author{Michel Balazard}
\date{}
\begin{document}
\maketitle

Le texte qui suit est une version détaillée de l'article {\og{\cyr \`Elementarnye zamechaniya o funkcii M{\"e}biusa}\fg}, paru aux {\cyr Trudy Matematicheskogo Instituta im. V.A. Steklova} {\bf 276} (2012), p.39-45 (voir \cite{Trudy} pour une traduction en langue anglaise).

\begin{flushright}
  {\it À la mémoire d'A. A. Karatsuba, pour le 75\up{e} anniversaire de sa naissance.}
\end{flushright}

\section{L'inversion de Möbius}

 En 1832, Möbius a défini la fonction arithmétique (qui porte maintenant son nom)
$$
\mu(n) =
\begin{cases}
  (-1)^r &\text{si $n$ est le produit de $r$ nombres premiers deux à deux distincts ;}\\
\quad 0 &\text{si $n$ est divisible par le carré d'un nombre premier}
\end{cases}
$$
(cf. \cite{009.0333cj}). La notation $\mu$ est due à Mertens (cf. \cite{06.0114.02})\footnote{Pour avoir une vue d'ensemble de l'histoire de la fonction de Möbius au XIX\up{e} siècle, il faut consulter l'indispensable ouvrage de Dickson \cite{MR0245499}, chapter XIX.}.

La propriété fondamentale de la fonction de Möbius, l'\emph{inversion de Möbius}, consiste en ceci. Soit $\Fcal$ l'ensemble des fonctions numériques définies sur $]0,\infty[$ et nulles sur $]0,1[$. On pose pour $\fhi \in \Fcal$
\begin{align*}
  S_{\un}\fhi(x) &=\sum_{n\soe 1}\fhi(x/n)\, ,\\
  S_{\mu}\fhi(x) &=\sum_{n\soe 1}\mu(n)\fhi(x/n) \quad (x>0),
\end{align*}
où les sommes ne portent en fait que sur les $n\ioe x$, par définition de $\Fcal$. Alors $S_{\un}$ et $S_{\mu}$ sont des applications linéaires de $\Fcal$ dans lui-même, inverses l'une de l'autre :
\begin{align*}
  \fhi(x) &=\sum_{n\soe 1}\mu(n)S_{\un}\fhi(x/n)\\
  \fhi(x) &=\sum_{n\soe 1}S_{\mu}\fhi(x/n) \quad (x>0).
\end{align*}

J'appelle \emph{transformation de Riemann} l'application $S_{\un}$ et \emph{transformation de Möbius} l'application $S_{\mu}$. Le nom de Riemann s'impose puisqu'en posant $h=1/x$ et $\omega(t)=\fhi(1/t)$ on voit que $x^{-1}S_{\un}\fhi(x)=h\sum_n\omega(nh)$ est une somme de Riemann de pas $h$ pour la fonction $\omega$ ; d'autre part, en calculant\footnote{On s'autorise ici une fonction $\fhi$ non nulle sur $]0,1[$, mais menant à une série $S_{\un}\fhi(x)$ absolument convergente. D'ailleurs, on peut considérer $S_{\un}$ et $S_{\mu}$ sur d'autres espaces que $\Fcal$, comme l'espace des fonctions nulles au voisinage de $0$, ou l'espace des fonctions qui sont $O(x^{\alpha})$ au voisinage de $0$ (avec $\alpha >1$).} $S_{\un}\fhi$ pour $\fhi(x)=x^s$ ($s>1$), on obtient 
$$
S_{\un}\fhi(x)=x^s\sum_n\frac{1}{n^s}=x^s\zeta(s),
$$
où $\zeta$ est la fonction de Riemann.

\section{Une identité de Meissel et une inégalité de Gram}\label{t28}

 L'application la plus simple de l'inversion de Möbius consiste à prendre\footnote{J'utilise la notation d'Iverson : $[A]=1$ si la proposition $A$ est vraie, et $[A]=0$ si elle est fausse.} $\fhi(x)=[x\soe 1]$, ce qui donne
\begin{align*}
S_{\un}\fhi(x)&=\sum_{n\ioe x}1\\
&=\lfloor x \rfloor\\
&=x-\{x\},  
\end{align*}
où $\lfloor x \rfloor$ et $\{x\}$ désignent respectivement la partie entière et la partie fractionnaire du nombre réel $x$. On obtient donc la formule de Meissel (1854, cf. \cite{048.1293cj}, p.303) :

\begin{equation}
  \label{t3}
\sum_{n\ioe x}\mu(n)\lfloor x/n\rfloor =1 \quad (x\soe 1).  
\end{equation}

Cette identité a été exploitée par Gram en 1884 (cf. \cite{16.0146.03}, p.197-198) pour démontrer le premier résultat non trivial sur la série $\sum \mu(n)/n$, à savoir que ses sommes partielles sont bornées :
$$
m(x)=\sum_{n\ioe x}\frac{\mu(n)}{n} \in [-1,1] \quad (x \soe 1).
$$

Rappelons la démonstration, très simple. Il suffit de traiter le cas de $x$ entier. On a alors
\begin{align}
  xm(x) &=\sum_{n\ioe x}\mu(n)\frac{x}{n}\notag\\
&=\sum_{n\ioe x}\mu(n)\lfloor x/n\rfloor +\sum_{n\ioe x}\mu(n)\{x/n\}\notag\\
&=1+\sum_{n\ioe x}\mu(n)\{x/n\}\label{t0}\\
&\in [2-x,x]\notag
\end{align}
puisque $|\mu(n)\{x/n\}|\ioe 1$ et $\{x/x\}=0$.

\section{Une identité de MacLeod faisant intervenir les polynômes de Bernoulli}

MacLeod a donné dans \cite{MR1301298} une généralisation de la formule de Meissel en termes des polynômes de Bernoulli $b_k$, définis par l'identité formelle
$$
\frac{z}{e^z-1}e^{zX}=\sum_{k\soe 0}b_k(X)\frac{z^k}{k!}.
$$
On a
\begin{align*}
  b_0(X) &= 1\\
b_1(X) &=X-\demi\\
b_2(X) &=X^2-X+\frac 16\\
b_3(X) &=X^3-\frac 32 X^2+\frac X2\\
b_4(X) &=X^4-2 X^3+X^2-\frac{1}{30}\\
b_5(X) &=X^5-\frac 52 X^4+\frac 53 X^3-\frac X6\\
b_6(X) &=X^6-3X^5 +\frac{5}{2}X^4-\frac{1}{2}X^2+\frac{1}{42}\\
b_7(X) &=X^7-\frac 72X^6 +\frac{7}{2}X^5-\frac{7}{6}X^3+\frac{1}{6}X\\
b_8(X) &=X^8 -4X^7+\frac{14}{3}X^6-\frac{7}{3}X^4+\frac{2}{3}X^2-\frac{1}{30},
\end{align*}
etc.

Plus généralement,
$$
b_k(X)=\sum_{j=0}^kB_j\binom{k}{j}X^{k-j},
$$
où $B_j=b_j(0)$ est le $j$\up{e} nombre de Bernoulli. Rappelons que $B_j=0$ si $j$ est impair et supérieur à 1 (voir par exemple \cite{MR1397498}, \S 6.5). 

La suite des $b_k$ vérifie de nombreuses identités, en particulier celle-ci :
\begin{equation}
  \label{t1}
b_k(X+1)-b_k(X)=kX^{k-1}.  
\end{equation}

Pour $k$ entier positif et $x$ réel positif, posons
$$
B_k(x)=b_k(\{x\}),
$$
et
$$
\fhi_k(x)=\frac{b_k(x)-B_k(x)}{x^{k-1}}.
$$

Les fonctions $\fhi_k$ appartiennent à $\Fcal$. On a
\begin{align*}
 \fhi_1(x)&=x-\{x\}=\lfloor x \rfloor\\
\fhi_2(x) &=x-1-\frac{\{x\}^2-\{x\}}{x}\\
\fhi_3(x) &=x-\frac 32 +\frac{1}{2x}-\frac{\{x\}^3-\frac 32\{x\}^2+\demi\{x\}}{x^2}\\ 
\fhi_4(x) &=x-2 +\frac{1}{x}-\frac{\{x\}^4-2\{x\}^3+\{x\}^2}{x^3}\\
\fhi_5(x) &=x-\frac 52 +\frac{5}{3x}-\frac{1}{6x^3}-\frac{\{x\}^5-\frac 52\{x\}^4+\frac 53 \{x\}^3-\frac 16\{x\}}{x^4}\\ 
\fhi_6(x) &=x-3 +\frac{5}{2x}-\frac{1}{2x^3}-\frac{\{x\}^6-3\{x\}^5+\frac 52\{x\}^4-\demi\{x\}^2}{x^5}\\ 
\fhi_7(x) &=x-\frac 72 +\frac{7}{2x}-\frac{7}{6x^3}+\frac{1}{6x^5}-\frac{\{x\}^7-\frac 72\{x\}^6+\frac 72\{x\}^5-\frac 76\{x\}^3+\frac 16\{x\}}{x^6}\\ 
\fhi_8(x) &=x-4 +\frac{14}{3x}-\frac{7}{3x^3}+\frac{2}{3x^5}-\frac{\{x\}^8-4\{x\}^7+\frac{14}{3}\{x\}^6-\frac{7}{3}\{x\}^4 +\frac 23\{x\}^2}{x^7},
\end{align*}
etc.

Plus généralement,
\begin{equation}
  \label{t25}
\fhi_k(x)=\sum_{j=0}^{k-1}B_j\binom{k}{j}x^{1-j} -\frac{\sum_{j=0}^{k-1}B_j\binom{k}{j}\{x\}^{k-j}}{x^{k-1}}.  
\end{equation}

La proposition suivante est due à MacLeod. La formulation et la démonstration que j'en donne ici sont différentes de celles de \cite{MR1301298}.

\begin{prop}[MacLeod, 1994]\label{t4}
Pour $k$ entier positif et $x\soe 1$, on a
$$
\sum_{n\ioe x}\mu(n)\fhi_k(x/n)=k(1-1/x)^{k-1}.
$$  
\end{prop}
\dem

Par inversion de Möbius, il s'agit de montrer que $\fhi_k(x)$ est la transformée de Riemann de la fonction $k(1-1/x)^{k-1}[x\soe 1]$.

Effectivement :
\begin{align*}
  \sum_{n\ioe x}k(1-n/x)^{k-1} &= x^{1-k}\sum_{n\ioe x}k(x-n)^{k-1}\\
&=x^{1-k}\sum_{n\ioe x}\big ( b_k(x-n+1)-b_k(x-n)\big)\expli{en vertu de \eqref{t1}}\\
&=x^{1-k}\big ( b_k(x-1+1)-b_k(x-\lfloor x \rfloor)\big)\expli{somme télescopique}\\
&=x^{1-k}\big ( b_k(x)-b_k(\{x \})\big)\\
&=\fhi_k(x).\fine
\end{align*}

\medskip

Terminons ce paragraphe avec quelques observations sur les fonctions $\fhi_k$. On a d'abord
\begin{align*}
 \fhi_k(x)&=\frac{b_k(x)-b_k(\{x\})}{x^{k-1}}\\
&=kx^{1-k}\sum_{n\ioe x} (x-n)^{k-1}, 
\end{align*}
où la deuxième égalité a été vue lors de la démonstration de la proposition \ref{t4}.

Rappelons que les \emph{moyennes intégrales de Cesàro}\footnote{Hardy, dans \cite{H}, \S 5.14, définit les moyennes intégrales de Cesàro d'une \emph{fonction de variable réelle} $a(t)$ par $A_0(x)=\int_0^xa(t)dt$ et $A_k(x)=\int_0^xA_{k-1}(t)dt$. Cette définition se généralise en remplaçant $a(t)dt$ par une \emph{mesure} borélienne quelconque (finie sur les segments), et les moyennes dont nous parlons sont alors celles de la mesure $\sum_nf(n)\delta_n$, classiquement associée à la fonction arithmétique $f$.} d'une fonction arithmétique $f$ sont définies par la formule
$$
F_k(x) =\sum_{n\ioe x}f(n) \frac{(x-n)^k}{k!} \quad (k\in \Nat,\, x>0).
$$
Ainsi, $F_0$ est la fonction sommatoire de $f$ et
$$
F_{k+1}(x)=\int_0^xF_k(t)dt.
$$

On voit donc que $x^{k-1}\fhi_k(x)/k!$ est la $(k-1)$\up{e} moyenne intégrale de Cesàro de la fonction arithmétique constante égale à $1$. En particulier, on a l'identité
$$
x^k\frac{\fhi_{k+1}(x)}{(k+1)!}=\int_0^xt^{k-1}\frac{\fhi_{k}(t)}{k!}dt,
$$
\cad 
$$
\fhi_{k+1}(x)=(k+1)x^{-k}\int_0^xt^{k-1}\fhi_{k}(t)dt.
$$

On a donc aussi
$$
\frac{d}{dx}\big (x^k\fhi_{k+1}(x)\big) = (k+1)x^{k-1}\fhi_{k}(x).
$$

Le comportement de la fonction $\fhi_k$ aux voisinages de $1$ et de l'infini est donné en première approximation par
$$
\fhi_k(x)\sim k(x-1)^{k-1} \quad (x\vers 1) \et \fhi_k(x)=x-k/2+O(1/x) \quad (x\vers \infty).
$$

Enfin, voici une identité qui nous sera utile dans les applications de \eqref{t20} :
\begin{equation}
  \label{t21}
  \int_1^{\infty}B_k(t)\frac{dt}{t^{k+1}}=1+\sum_{j=1}^k\frac{B_j}{j}-\gamma, 
\end{equation}
où $\gamma$ est la constante d'Euler. Cette identité se démontre par intégrations par parties successives à partir du cas $k=1$ (voir aussi \cite{MR1732751}, Theorem 2.2).

\section{Une identité générale}

 Soit $f$ une fonction arithmétique, $F$ sa fonction sommatoire, et $\fhi :[1,\infty[ \vers \Real$ une fonction absolument continue sur tout segment, s'annulant au point $1$. Dans la suite, je ferai plusieurs fois usage de l'identité
\begin{align}
  S_f\fhi(x)= \sum_{n\ioe x}f(n)\fhi(x/n) &=\sum_{n\ioe x}f(n)\int_1^{x/n}\fhi'(t)dt\notag\\
&=\int_1^{x}F(x/t)\fhi'(t)dt,\label{t24}
\end{align}
où l'on a interverti sommation et intégration pour obtenir la dernière égalité.

En particulier, en supposant
$$
\alpha=\limsup_{x\vers \infty}x^{-1}|F(x)| <\infty,
$$
on a
$$
\limsup_{x\vers \infty}x^{-1}|S_f\fhi(x)|\ioe \alpha\int_1^{\infty}|\fhi'(t)|\frac{dt}{t},
$$ 
pourvu que cette dernière intégrale soit convergente. Si en outre $\fhi'\soe 0$, cela donne
\begin{equation}
  \label{t20}
 \limsup_{x\vers \infty}x^{-1}|S_f\fhi(x)|\ioe \alpha\int_1^{\infty}\fhi'(t)\frac{dt}{t}=\alpha\int_1^{\infty}\fhi(t)\frac{dt}{t^2}. 
\end{equation}

\section{Le cas $k=2$ de l'identité de Macleod}\label{t30}

 Le cas $k=1$ de la proposition \ref{t4} est la formule de Meissel \eqref{t3}. Examinons maintenant le cas $k=2$ :
\begin{equation}
  \label{t5}
  \sum_{n\ioe x} \mu(n)\Big (\frac xn-1-\frac{\{x/n\}^2-\{x/n\}}{x/n}\Big )=2-\frac 2x \quad (x\soe 1).
\end{equation}

Notons suivant l'usage
$$
M(x)=\sum_{n\ioe x}\mu(n) 
$$
la fonction sommatoire ordinaire de la fonction de Möbius (la fonction $m(x)$ est sa fonction sommatoire \emph{logarithmique}). Posons également pour $x>0$ :
\begin{align}
  m_1(x) &=\int_1^x M(t)\frac{dt}{t^2}\notag\\
&=\frac 1x \int_1^xm(t)dt\label{t22}\\
&=m(x)-\frac{M(x)}{x}\label{t23}\\
&=\sum_{n\ioe x}\mu(n) \Big (\frac 1n -\frac 1x\Big),\notag
\end{align}
ces identités étant d'ailleurs valables généralement, en y remplaçant $\mu$ par n'importe quelle autre fonction.

En appliquant \eqref{t24}, on peut récrire \eqref{t5} sous la forme

\begin{align}
  xm_1(x) &= 2-\frac 2x +\sum_{n\ioe x} \mu(n)\frac{\{x/n\}^2-\{x/n\}}{x/n}\notag\\
&=2-\frac 2x +\int_1^{x} M(x/t)\frac{(2\{t\}-1)t+\{t\}-\{t\}^2}{t^2}dt.\label{t6}
\end{align}

\begin{prop}\label{t7}
Pour $t\soe 0$, on a
$$
|(2\{t\}-1)t+\{t\}-\{t\}^2|\ioe t.
$$
\end{prop}
\dem

On a d'abord
$$
(2\{t\}-1)t+\{t\}-\{t\}^2\soe -t.
$$
Ensuite,
\begin{align*}
  (2\{t\}-1)t+\{t\}-\{t\}^2 -t &=(\{t\}-1)(2t-\{t\})\\
&\ioe 0.\fine
\end{align*}

La conjonction de \eqref{t6} et de la proposition \ref{t7} fournit le résultat suivant.
\begin{prop}\label{t8}
Pour $x\soe 1$, on a
$$
|m_1(x)|\ioe \frac 1x \int_1^x|M(t)|\frac {dt}{t} +\frac 2x -\frac{2}{x^2}.
$$ 
\end{prop}
\dem

On a
\begin{align*}
  \Big |\int_1^{x} M(x/t)\frac{(2\{t\}-1)t+\{t\}-\{t\}^2}{t^2}dt\Big | &\ioe \int_1^{x} |M(x/t)|\frac {dt}{t}\expli{d'aprés la proposition \ref{t7}}\\
&= \int_1^x|M(t)|\frac {dt}{t}.\fine
\end{align*}

En particulier, en posant 
$$
\alpha=\limsup_{x\vers \infty} |M(x)/x|,
$$
on obtient
$$
\limsup_{x\vers \infty} |m_1(x)|\ioe \alpha
$$
et donc
\begin{equation}
  \label{t9}
  \limsup_{x\vers \infty} |m(x)|\ioe 2\alpha.
\end{equation}

\section{Le lemme d'Axer}\label{t38}

 La première démonstration de l'implication
\begin{equation}
  \label{t10}
  M(x)=o(x)\implique m(x)=o(1) \quad \quad (x\vers \infty) 
\end{equation}
par les méthodes de l'analyse réelle\footnote{La démonstration de $m(x)=o(1)$ par les méthodes de l'analyse complexe fut obtenue dés 1897 par von Mangoldt (cf. \cite{28.0180.01}) et ce résultat implique immédiatement que $M(x)=o(x)$, d'aprés un théoréme général de Kronecker sur les séries convergentes (cf. \cite{18.0212.04}), ou bien par l'égalité de \eqref{t22} et \eqref{t23}.} fut obtenue par Axer en 1910 (cf. \cite{41.0240.02}) ; c'est encore de nos jours celle qui est présentée dans les manuels de théorie analytique des nombres. Le raisonnement présenté au \S 4 constitue une autre démonstration de ce résultat ; il présente en outre un avantage quantitatif sur celui d'Axer, comme nous allons le voir maintenant. 

Le lemme fondamental utilisé par Axer est le suivant.
\begin{prop}[Axer, 1910]
Soit $f$ une fonction arithmétique,
$$
F(x)=\sum_{n\ioe x}f(n) \quad \text{et}\quad G(x)=\sum_{n\ioe x}|f(n)|.
$$
On suppose que
$$
F(x)=o(x) \quad (x\vers \infty)\quad \text{et}\quad G(x)=O(x) \quad (x\soe 1).
$$
Alors
$$
\sum_{n\ioe x}f(n)\{x/n\}=o(x) \quad (x\vers \infty).
$$  
\end{prop}

La conjonction de la formule de Meissel sous la forme \eqref{t0} et du lemme d'Axer fournit bien l'implication \eqref{t10}.

La lecture de la démonstration d'Axer, ou des démonstrations ultérieures de Landau, en révèle des variantes quantitatives, comme la suivante (cf. \cite{41.0241.01}, p.132-134).
\begin{prop}[Landau, 1910]\label{t39}
Soit $f$ une fonction arithmétique,
$$
F(x)=\sum_{n\ioe x}f(n) \quad \text{et}\quad G(x)=\sum_{n\ioe x}|f(n)|.
$$
On suppose que
$$
\limsup_{x\vers \infty} |F(x)/x| \ioe \alpha \ioe 1/2  \quad \text{et}\quad |G(x)|\ioe Ax \quad (x\soe 1).
$$
Alors
$$
\limsup_{x\vers \infty}\big |x^{-1}\sum_{n\ioe x}f(n)\{x/n\}\big |\ioe \alpha \big (A+5+2\log (1/\alpha)\big ).
$$  
\end{prop}

Nous verrons au \S \ref{t41} que la présence du terme $\alpha\log (1/\alpha)$ est inévitable à ce degré de généralité. Pour $0<\alpha\ioe 1/2$, l'application de cette version quantitative du lemme d'Axer donne 
\begin{equation*}
    \limsup_{x\vers \infty} |m(x)|\ll \alpha\log (1/\alpha)
\end{equation*}
au lieu de \eqref{t9}.

\section{Application des cas $k=3$ et $k=4$ de l'identité de MacLeod}

 En passant maintenant à $k=3$, la proposition \ref{t4} nous donne
$$
\sum_{n\ioe x}\mu(n)\Big (\frac xn -\frac 32 +\frac{1}{2x/n}- \frac{\{x/n\}^3-\frac 32 \{x/n\}^2+\demi\{x/n\}}{(x/n)^2}\Big)=3(1-1/x)^2\quad (x\soe 1).
$$

On peut à partir de cette identité démontrer la proposition \ref{t14} ci-dessous (avec $3$ au lieu de $8/3$). Je préfère cependant utiliser la fonction suivante, découverte après quelques essais\footnote{Pour les expériences numériques effectuées lors de la préparation de cet article, j'ai utilisé le logiciel Maple, version 9.5.}.

 \begin{equation*}
  \frac 43 \fhi_3(x)- \frac 13 \fhi_4(x) =x-1-\eps_1(x),
\end{equation*}
où
$$
\eps_1(x)=\frac 13-\frac{1}{3x}+\frac 43 \frac{\{x\}^3-\frac 32\{x\}^2+\demi\{x\}}{x^2} -\frac 13\frac{\{x\}^4-2\{x\}^3+\{x\}^2}{x^3}.
$$

La fonction $\eps_1$ vérifie une identité remarquable.

\begin{prop}\label{t17}
Pour $x>0$, on a $\eps'_1(x)=\big(1-\fhi'_2(x)\big )^2$.  
\end{prop}
\dem

En effet,
\begin{align*}
\eps'_1(x)&=\frac{1}{3x^2}+4\frac{\{x\}^2-\{x\}+\frac 16}{x^2} -\frac 83\frac{\{x\}^3-\frac 32\{x\}^2+\demi\{x\}}{x^3} -\frac 43 \frac{\{x\}^3-\frac 32\{x\}^2+\demi\{x\}}{x^3}+\frac{\{x\}^4-2\{x\}^3+\{x\}^2}{x^4}\\
&= \frac{(2\{x\}-1)^2}{x^2}-2\frac{\{x\}(\{x\}-1)(2\{x\}-1)}{x^3}+\frac{\{x\}^2(\{x\}-1)^2}{x^4}\\
&= \Big (\frac{(2\{x\}-1)x-\{x\}(\{x\}-1)}{x^2}\Big )^2\\
&=\big(1-\fhi'_2(x)\big )^2.\fine
\end{align*}

La proposition \ref{t17} peut se récrire sous la forme
$$
3-4\fhi'_3+\fhi'_4=3(1-\fhi'_2)^2.
$$
Observons que la fonction $\fhi'_2$ a une discontinuité en chaque entier, mais que $(1-\fhi'_2)^2$ est continue. 

\begin{prop}\label{t14}
Pour $x\soe 1$, on a
$$
|m_1(x)|\ioe \frac{1}{x^2}\int_1^x|M(t)|dt+8/3x.
$$
\end{prop}
\dem

On a
\begin{align*}
xm_1(x)&=\sum_{n\ioe x}\mu(n)\Big(\frac xn -1\Big)\\
&=\sum_{n\ioe x}\mu(n)\Big(\frac 43\fhi_3(x/n)-\frac 13\fhi_4(x/n)\Big)+\sum_{n\ioe x}\mu(n)\eps_1(x/n)\\
&=\frac 43 \cdot 3(1-1/x)^2-\frac 13 \cdot 4(1-1/x)^3+\int_1^xM(x/t)\eps_1'(t)dt,
\end{align*}
d'après la proposition \ref{t4} et l'identité \eqref{t24}.

On montre d'une part facilement que 
$$
0\ioe \frac 43\cdot3(1-1/x)^2-\frac 13\cdot4(1-1/x)^3 \ioe \frac 83 \quad (x\soe 1).
$$

D'autre part, les propositions \ref{t7} et \ref{t17} entraînent que $0\ioe \eps_1'(t) \ioe t^{-2}$. Le résultat en découle.\fin

\medskip

En particulier, on déduit de la proposition \ref{t14} que
\begin{equation}
  \label{t15}
\limsup_{x\vers \infty}|m_1(x)|\ioe \demi \limsup_{x\vers \infty}|M(x)/x|.\footnote{Bien entendu, on sait depuis 1897 que ces deux limites sont nulles, mais la démonstration que nous venons de donner de cette inégalité ne dépend pas de ce fait.}
\end{equation}

Cela étant, d'aprés la proposition \ref{t17}, $\eps'_1(x)$ est positif ou nul. Les relations \eqref{t20} et \eqref{t21} montrent alors qu'on peut remplacer la constante $\demi$ de \eqref{t15} par
\begin{align*}
\int_1^{\infty}\eps_1(t)\frac{dt}{t^2}&=\int_1^{\infty}\Bigl ( \frac{1}{3t^2}-\frac{1}{3t^3}+\frac 43 \frac{B_3(t)}{t^4}-\frac 13 \frac{B_4(t)}{t^5}\Bigr) )dt\\
&= \frac{271}{360}-\gamma=0,1755\dots
\end{align*}

\section{Généralisation des propositions \ref{t8} et \ref{t14}}\label{t44}

Je donne maintenant une proposition générale, évidente pour $k=1$, et dont les propositions \ref{t8} et \ref{t14} explicitent les cas particuliers $k=2$ et $k=3$.
\begin{prop}\label{t18}
  Pour tout nombre entier $k$, il existe deux constantes positives $C_k$ et $D_k$ telles que
$$
|m_1(x)|\ioe C_kx^{1-k}\int_1^x|M(t)|t^{k-3}dt +D_k/x \quad (x>0).
$$
\end{prop}
\dem

Comme $m_1(x)=0$ si $0<x<1$, on peut supposer $x\soe 1$. Observons d'abord que le résultat pour une valeur particulière de $k$ entraîne le résultat pour toutes les valeurs inférieures (avec les mêmes constantes $C_k$ et $D_k$) puisque pour $j\ioe k$,
$$
x^{-k}t^{k-2}\ioe x^{-j}t^{j-2} \quad (1\ioe t\ioe x).
$$ 
On peut donc supposer $k$ entier, impair et supérieur ou égal à $5$, disons. Avec un léger changement de notation, nous allons donc montrer que, pour $k$ entier supérieur ou égal à $2$,
\begin{equation}
  \label{t19}
|m_1(x)|\ll_k x^{-2k}\int_1^x|M(t)|t^{2k-2}dt +x^{-1} \quad (x\soe 1).  
\end{equation}

Pour cela, nous considérons la fonction
$$
\fhi=\lambda_1\fhi_{2k+1}+\lambda_2\fhi_{2k+2}+\cdots+\lambda_k\fhi_{3k},
$$
avec des coefficients $\lambda_1,\dots,\lambda_k$ convenables.

On a
\begin{align*}
  \fhi(x) &=\sum_{\ell=1}^k \lambda_{\ell}\fhi_{2k+\ell}(x)\\
&= \sum_{\ell=1}^k \lambda_{\ell}\Bigl ( x-\demi\binom{2k+\ell}{1}+\sum_{i=1}^{k-1}B_{2i}\binom{2k+\ell}{2i}x^{1-2i}\Bigr)+\\
&\quad +\sum_{\ell=1}^k \lambda_{\ell}\sum_{i=k}^{k+\lfloor(\ell-1)/2\rfloor}B_{2i}\binom{2k+\ell}{2i}x^{1-2i}-\sum_{\ell=1}^k \lambda_{\ell}\frac{B_{2k+\ell}(x)-B_{2k+\ell}}{x^{2k+\ell-1}},
\end{align*}
où l'on a séparé les termes correspondant à $j=0$, $j=1$, $j=2i$ (où $1\ioe i\ioe k-1$), et $j=2i$ (où $i\soe k$) dans la première somme de \eqref{t25} (avec $k$ remplacé par $2k+\ell$).

Le terme
$$
\psi(x)=\sum_{\ell=1}^k \lambda_{\ell}\sum_{i=k}^{k+\lfloor(\ell-1)/2\rfloor}B_{2i}\binom{2k+\ell}{2i}x^{1-2i}-\sum_{\ell=1}^k \lambda_{\ell}\frac{B_{2k+\ell}(x)-B_{2k+\ell}}{x^{2k+\ell-1}}$$
vérifie $\psi'(x)\ll_kx^{-2k}$ (où la constante implicite dépend de $k$ et des $\lambda_{\ell}$). On a ensuite
$$
\fhi(x)=x\sum_{\ell=1}^k \lambda_{\ell} -\sum_{\ell=1}^k (k+\ell/2)\lambda_{\ell}+\sum_{i=1}^{k-1}B_{2i}\frac{\sum_{\ell=1}^k \binom{2k+\ell}{2i}\lambda_{\ell}}{x^{2i-1}}+\psi(x),
$$
et on cherche les $\lambda_{\ell}$ de sorte que
\begin{align}
  \sum_{\ell=1}^k \lambda_{\ell} &=1\label{t26}\\
\sum_{\ell=1}^k \binom{2k+\ell}{2i}\lambda_{\ell}&=0 \quad (1\ioe i \ioe k-1).\label{t27} 
\end{align}
Le déterminant de ce système de $k$ équations linéaires à $k$ inconnues est $\Delta_k(2k)$, où l'on a posé 
$$
\Delta_k(x)=\det \big (\binom{x+j}{2(i-1)}\big)_{1\ioe i,j\ioe k}.
$$

On a
$$
\binom{x+j}{2(i-1)}-\binom{x+j-1}{2(i-1)}=\binom{x+j-1}{2i-3}=\frac{x+j-1}{2i-3}\binom{x+j-2}{2(i-2)} \quad (2\ioe i,j\ioe k),
$$
et cette différence est nulle pour $i=1$. En retranchant la $(j-1)$\up{e} colonne de $\Delta_k(x)$ à sa $j$\up{e} pour $k\soe j\soe 2$, on voit donc que
$$
\Delta_k(x)= \frac{\prod_{j=2}^k (x+j-1)}{\prod_{i=2}^k(2i-3)}\Delta_{k-1}(x-1).
$$

Comme $\Delta_1(x)=1$, on en déduit par récurrence sur $k$ que $\Delta_k(x)\not= 0$ pour $x\soe k$, disons. En particulier $\Delta_k(2k)\not =0$.

On peut donc choisir les $\lambda_{\ell}$ tels que les relations \eqref{t26} et \eqref{t27} soient vérifiées. Nous aurons alors
$$
\fhi(x)=x-c_k +\psi(x),
$$
où $c_k=\sum_{\ell=1}^k (k+\ell/2)\lambda_{\ell}$. Comme $\fhi(1)=0$, on en déduit que
$$
x-1=\fhi(x)-\int_1^x\psi'(t)dt,
$$
donc 
\begin{align*}
  xm_1(x) &=\sum_{n\ioe x}\mu(n)\fhi(x/n)-\sum_{n\ioe x}\mu(n)\int_1^{x/n}\psi'(t)dt\\
&=\sum_{\ell=1}^k \lambda_{\ell}(2k+\ell)(1-1/x)^{2k+\ell-1}-\int_1^xM(x/t)\psi'(t)dt\\
&\ll_k 1+\int_1^x|M(x/t)|t^{-2k}dt\\
&=1+x^{1-2k}\int_1^x|M(t)|t^{2k-2}dt,
\end{align*}
ce qui démontre \eqref{t19}.\fin

\medskip

On peut voir que les constantes $C_k$ et $D_k$ fournies par la démonstration de la proposition \ref{t18} sont $O(6^{k^2})$. D'autre part, quelques essais me conduisent à conjecturer que les valeurs suivantes sont admissibles dans la proposition  \ref{t18}: 

\begin{center}
\begin{tabular}{|c||c|c|c|c|c|c|}
\hline
$k$&4&5&6&7&8\\
\hline
$C_k$ & 1,1&1,5&2,6&6,3&13,8\\
\hline
$D_k$ &2,1&2,5&2,8&3,1&3,2\\
\hline
\end{tabular}  
\end{center}

\section{Sur un encadrement dû à von Mangoldt}

Comme von Mangoldt l'a montré dans \cite{28.0180.01} (Hülfssatz 2, p.839), on peut adapter l'argumentation de Gram  à une identité de départ autre que celle de Meissel, mais toujours obtenue par inversion de Möbius. Au lieu de $\fhi(x)=[x\soe 1]$ comme au \S \ref{t28}, prenons cette fois $\fhi(x)=x[x\soe 1]$. Nous aurons
\begin{align*}
  S_{\un}\fhi(x)&=\sum_{n\ioe x}x/n\\
&=xH(x),
\end{align*}
avec
\begin{align*}
  H(x) &=\sum_{n\ioe x}1/n\\
&=\log x +\gamma +\frac{1/2-\{x\}}{x}+\eps_2(x),
\end{align*}
et
$$
\eps_2(x) =\int_x^{\infty}(\{t\}-1/2)\frac{dt}{t^2}
$$
d'après la formule sommatoire d'Euler et Maclaurin (cf. \cite{MR2382467}, proposition 1).

Par inversion de Möbius, on en déduit
$$
\sum_{n\ioe x}\mu(n)(x/n)H(x/n)=x \quad (x\soe 1),
$$
autrement dit :
$$
\sum_{n\ioe x}\frac{\mu(n)}{n}\big(\log (x/n)+\gamma+\frac{1/2-\{x/n\}}{x/n}+\eps_2(x/n)\big)=1.
$$

Par conséquent, on obtient
$$
\sum_{n\ioe x}\frac{\mu(n)}{n}\log (x/n)=\int_1^xm(t)\frac{dt}{t}=1-\gamma m(x)-\sum_{n\ioe x}\frac{\mu(n)}{n}\Big(\frac{1/2-\{x/n\}}{x/n}+\eps_2(x/n)\Big).
$$

En utilisant les inégalités $|m(x)|\ioe 1$, $|\eps_2(x)|\ioe 1/2x$ et $|\mu(n)|\ioe 1$, on obtient que
$$
\big |\gamma m(x)+\sum_{n\ioe x}\frac{\mu(n)}{n}\Big(\frac{1/2-\{x/n\}}{x/n}+\eps_2(x/n)\Big )\big | \ioe \gamma +1,
$$
d'où\footnote{Von Mangoldt donne $-3-\gamma$ et $3+\gamma$ comme bornes de cet encadrement. On peut remplacer les bornes $-\gamma$ et $2+\gamma$ de \eqref{t29} par les valeurs optimales $0$ et $\sum_{n\ioe 30}\frac{\mu(n)}{n}\log (30/n)=1,00302\dots$}
  \begin{equation}
    \label{t29}
-\gamma \ioe \sum_{n\ioe x}\frac{\mu(n)}{n}\log (x/n)\ioe 2+\gamma.    
  \end{equation}

Ainsi non seulement $m(t)$ est bornée, mais c'est aussi le cas de l'intégrale $\int_1^xm(t)dt/t$. 

\smallskip

Maintenant, de même que la proposition \ref{t8} du \S\ref{t30} donnait une majoration non triviale de $|m(t)|$ directement en termes de $|M(t)|$, de même on peut donner une majoration non triviale de $|-1+\int_1^xm(t)dt/t|$ en termes de $|M(t)|$ : c'est l'objet de la proposition \ref{t33} ci-dessous. La méthode est la même : au lieu de considérer comme au \S\ref{t30}
$$
\fhi_2(x)=\frac 2x \int_0^x(\sum_{n\ioe t}1) dt,
$$
on considère maintenant
\begin{align*}
\beta_2(x)&=\frac 2x \int_0^x(\sum_{n\ioe t}t/n) dt\\
&=\sum_{n\ioe x}(x/n-n/x)\\
&=x(\log x+\gamma-1/2)+x\eps_2(x)-\frac{\{x\}^2-\{x\}}{2x}.  
\end{align*}

Par inversion de Möbius on obtient
\begin{equation}
  \label{t31}
 \sum_{n\ioe x}\mu(n)\big ((x/n)\log (x/n)+(\gamma-1/2)x/n+(x/n)\eps_2(x/n)-\frac{\{x/n\}^2-\{x/n\}}{2x/n}\big)=x-1/x \quad (x\soe 1). 
\end{equation}
La combinaison linéaire \eqref{t31}$+(1/2-\gamma)\times$\eqref{t5} fournit alors
\begin{equation}
  \label{t32}
  \sum_{n\ioe x}\mu(n)\big ((x/n)\log (x/n)-\eps_3(x/n)\big) =x+1-2\gamma+(2\gamma-2)/x,
\end{equation}
où
$$
\eps_3(x)=\gamma-1/2+x\eps_2(x)+(\gamma-1)\frac{\{x\}^2-\{x\}}{x}.
$$
Observons que $\eps_3(1)=0$ et que la fonction $\eps_3$ est continûment dérivable sauf aux points entiers, avec
\begin{align}
\eps_3'(x) &=x\eps_2'(x)+\eps_2(x)+(\gamma-1)\frac{(2\{x\}-1)x-\{x\}^2+\{x\}}{x^2}\notag\\
&=H(x)-\log x -\gamma +(\gamma-1)\frac{(2\{x\}-1)x-\{x\}^2+\{x\}}{x^2}\quad (x\soe 1, \, x \not \in \Nat).\label{t35}  
\end{align}

\begin{prop}\label{t34}
 On a
$$
|x\eps_3'(x)| \ioe 1 \quad (x\soe 1, \, x \not \in \Nat).
$$ 
\end{prop}
\dem

Pour $x\soe 1, \, x \not \in \Nat$, on a d'aprés \eqref{t35}
$$
x\eps_3'(x)=xH(x)-x\log x -\gamma x +(2\gamma-2)(\{x\}-1/2)-(\gamma-1)\frac{\{x\}^2-\{x\}}{x},$$
donc
\begin{align*}
\frac{d}{dx}\big (x\eps_3'(x)\big) &=H(x)-\log x-1-\gamma+(2\gamma-2)- (\gamma-1)\frac{(2\{x\}-1)x-\{x\}^2+\{x\}}{x^2}\\
&\ioe  -\gamma +2\gamma-2+1-\gamma\expli{car $H(x)\ioe \log x +1$ et $(2\{x\}-1)x-\{x\}^2+\{x\}\ioe x \ioe x^2$}\\
&<0.
\end{align*}

Il suffit donc de vérifier que pour $n\in \Nat^*$ on a
\begin{align}
 n\eps_3'(n)  &\ioe 1\label{t36}\\
(n+1-0)\eps_3'(n+1-0)  &\soe -1\label{t37}.
\end{align}

L'inégalité \eqref{t36} se récrit 
$$
H(n)\ioe \log n+\gamma+\frac{\gamma}{n},
$$
ce qui résulte de 
\begin{align*}
H(n)- \log n-\gamma &= \frac{1}{2n}+\int_n^{\infty}(\{t\}-1/2)\frac{dt}{t^2}\\
&<\frac{1}{2n}. 
\end{align*}

Enfin l'inégalité \eqref{t37} se récrit
$$
H(n)\soe \log(n+1)+\frac{\gamma n}{n+1},
$$
soit encore
$$
\eps_2(n)\soe \log(1+1/n)-\frac{\gamma }{n+1}-\frac{1}{2n}. 
$$

Cela résulte des inégalités suivantes :

$\bullet$ $\eps_2(n)>-1/12n^2$ (cf. \cite{MR1397498}, formule (6.66)) ;

$\bullet$ $\log(1+1/n)\ioe 1/n$ ;

$\bullet$ $(\gamma-\demi)n^2-7n/12-1/12 \soe 0$ pour $n\soe 8$, 

et d'une comparaison directe pour $1\ioe n\ioe 7$.\fin

\begin{prop}\label{t33}
Pour $x\soe 1$, on a
$$
\Big \lvert -1+\int_1^xm(t)\frac{dt}{t}\Big \rvert\ioe \frac 1x + \frac 1x\int_1^x|M(t)|\frac {dt}{t}.
$$  
\end{prop}
\dem

L'identité \eqref{t32} se récrit
\begin{align*}
\int_1^xm(t)\frac{dt}{t}&= \sum_{n\ioe x}\frac{\mu(n)}{n}\log (x/n)\\
&=1+\frac{1-2\gamma}{x}+\frac{2\gamma-2}{x^2}+\frac 1x\sum_{n\ioe x}\mu(n)\eps_3(x/n).  
\end{align*}

On a d'abord
$$
-\frac 1x\ioe \frac{1-2\gamma}{x}+\frac{2\gamma-2}{x^2} \ioe 0 \quad (x\soe 1),
$$
puis 
\begin{align*}
\Bigl \lvert \sum_{n\ioe x}\mu(n)\eps_3(x/n)\Big \rvert &= \Bigl \lvert \int_1^xM(x/t)\eps'_3(t)dt\Big \rvert\\
&\ioe \int_1^x|M(x/t)|\frac{dt}{t} \expli{d'aprés la proposition \ref{t34}}\\
&=\int_1^x|M(t)|\frac {dt}{t}.\fine  
\end{align*}

\medskip

On pourrait généraliser la proposition \ref{t33} en une analogue de la proposition \ref{t18}, et d'autre part considérer la moyenne logarithmique de Riesz d'ordre $k$ de la fonction $\mu(n)/n$ :
$$
\Rcal_k(x)=\frac{1}{k!}\sum_{n\ioe x}\frac{\mu(n)}{n}\log^k(x/n),
$$
dont le terme principal est $P_k(\log x)/k!$, où $P_k$ est le $k$\up{e} polynôme d'Appell de la fonction $1/\zeta(s)$ en $s=1$ (voir \cite{MR2382467}, \S 2). 

\section{Un exemple à propos du lemme d'Axer}\label{t41}

Dans la version de ces remarques soumise pour publication à l'Institut Steklov, je posais une question concernant le lemme d'Axer. Afin de la formuler définissons, pour $0\ioe \alpha\ioe 1/2$, la quantité $c(\alpha)$ comme la borne inférieure des $\beta>0$ tels que l'énoncé suivant soit vrai.
\begin{quote}
  \textit{Soit $f$ une fonction arithmétique,
$$
F(x)=\sum_{n\ioe x}f(n) \quad \text{et}\quad G(x)=\sum_{n\ioe x}|f(n)|.
$$
On suppose que
$$
\limsup_{x\vers \infty} |F(x)/x| \ioe \alpha   \quad \text{et}\quad \limsup_{x\vers \infty} G(x)/x \ioe 1.
$$
Alors
$$
\limsup_{x\vers \infty}\big |x^{-1}\sum_{n\ioe x}f(n)\{x/n\}\big |\ioe \beta.
$$  
}
\end{quote}

La proposition \ref{t39} entraîne que $c(\alpha)\ioe \alpha\big (6+2\log(1/\alpha)\big)$.
\begin{ques}\label{t49}
  A-t-on
$$
c(\alpha) =o\big(\alpha\log(1/\alpha)\big) \quad (\alpha \vers 0)?
$$
\end{ques}

L'arbitre anonyme sollicité par l'Institut Steklov a montré que la réponse à cette question est négative grâce à la construction suivante, que je présente en décomposant son argumentation sous la forme de plusieurs propositions.

\begin{prop}\label{t46}
Soit $\alpha >0$ tel que $\alpha^{-1}$ soit un entier supérieur ou égal à $2$. Alors on a
\begin{align*}
1) \quad &\alpha^{-4-2j} > \alpha^{-4-2j}/2+1 \quad (j\soe 0)\\
2) \quad &\lfloor \alpha^{-4-2(j+1)}/k\rfloor > \lfloor \alpha^{-4-2j}/k'\rfloor +1 \quad (j\soe 0 \, ; \, 2\ioe k,k'\ioe \alpha^{-1})\\
3) \quad &\lfloor \alpha^{-4-2j}/k\rfloor > \lfloor \alpha^{-4-2j}/(k+1)\rfloor +1 \quad (j\soe 0 \, ; \, 2\ioe k\ioe \alpha^{-1}-1).
\end{align*}
\end{prop}
\dem

Pour $1)$ on a en fait $\alpha^{-4-2j}/2 \soe 8>1$. Pour $2)$ on a
$$
\lfloor \alpha^{-4-2(j+1)}/k\rfloor \soe \lfloor \alpha^{-4-2(j+1)}/(\alpha^{-1})\rfloor=\alpha^{-5-2j}>\alpha^{-4-2j}> \alpha^{-4-2j}/2+1 \soe \lfloor \alpha^{-4-2j}/k'\rfloor +1.
$$

Enfin, pour $3)$, il suffit de voir que
$$
\alpha^{-4-2j}/\big(k(k+1)\big )\soe 2 \quad (j\soe 0,\; 2\ioe k\ioe \alpha^{-1}-1),
$$
et cela résulte de 
$$
\alpha^{-4}/\big(\alpha^{-1}(\alpha^{-1}-1)\big) \soe 2,
$$
pour $\alpha^{-1}\soe 2$.\fin

\begin{prop}\label{t42}
Soit $\alpha >0$ tel que $\alpha^{-1}$ soit un entier supérieur ou égal à $2$. Alors l'application
\begin{align*}
\{0,1\}\times\Nat\times\{k\in \Nat, \, 2\ioe k\ioe \alpha^{-1}\} &\vers \Nat^*\\
(h,j,k)&\mapsto\lfloor \alpha^{-4-2j}/k\rfloor+h
\end{align*}
est injective.
\end{prop}
\dem

La proposition \ref{t46} montre que cette application est strictement croissante si l'on munit l'ensemble de départ de l'ordre lexicographique sur les triplets $(j,k^{-1},h)$.\fin

\smallskip

La proposition \ref{t42} nous permet de définir, pour tout $\alpha$ inverse d'entier supérieur ou égal à $2$, une fonction arithmétique $f_{\alpha}$ par la formule
$$
f_{\alpha}(n)=
\begin{cases}
(-1)^h\alpha \lfloor \alpha^{-4-2j}/k\rfloor &\text{si $n=\lfloor \alpha^{-4-2j}/k\rfloor+h$} \quad(j\soe 0\; ; \; h=0,1\; ;\;2\ioe k\ioe \alpha^{-1})\\
0&\text{sinon.}
\end{cases}
$$

En désignant par $F_{\alpha}$ la fonction sommatoire de $f_{\alpha}$, on a pour $N\in \Nat^*$
\begin{equation}\label{t45}
F_{\alpha}(N) =
\begin{cases}
\alpha N &\text{si $N$ est de la forme }\lfloor \alpha^{-4-2j}/k\rfloor\quad(j\soe 0\; ;\;2\ioe k\ioe \alpha^{-1})\\
0&\text{sinon.}
\end{cases}
\end{equation}
En particulier on a
$$
\limsup |F_{\alpha}(x)/x|=\alpha.
$$

Posons ensuite
$$
G_{\alpha}(x)=\sum_{n\ioe x}|f_{\alpha}(n)| \quad (x>0).
$$
\begin{prop}\label{t43}
Soit $\alpha >0$ tel que $\alpha^{-1}$ soit un entier supérieur ou égal à $2$. Pour tout $x>0$, on a
$$
G_{\alpha}(x)\ioe (x+1)(2/e+3\alpha).
$$
\end{prop}
\dem

Au vu de la définition de la fonction arithmétique $f_{\alpha}$, il suffit de considérer le cas de $x=N=\lfloor \alpha^{-4-2J}/K\rfloor +h_0$, où $J\in \Nat$, $h_0\in\{0,1\}$ et $2\ioe K\ioe \alpha^{-1}$. 

On a
$$
G_{\alpha}(N)=\sum_{j=0}^{J-1}\sum_{k=2}^{\alpha^{-1}}2\alpha \lfloor \alpha^{-4-2j}/k\rfloor+\sum_{k=K+1}^{\alpha^{-1}}2\alpha \lfloor \alpha^{-4-2J}/k\rfloor +\alpha \lfloor \alpha^{-4-2J}/K\rfloor(h_0+1).
$$

Le dernier terme est $\ioe 2\alpha N$. Pour le deuxième, on a
\begin{align*}
\sum_{k=K+1}^{\alpha^{-1}}2\alpha \lfloor \alpha^{-4-2J}/k\rfloor &\ioe 2\alpha^{-3-2J}\sum_{k=K+1}^{\alpha^{-1}}1/k\\
&\ioe 2\alpha^{-3-2J}\int_K^{\alpha^{-1}}\frac{dt}{t}\\
&= 2\alpha^{-3-2J}\log (1/K\alpha)\\
&=2(\alpha^{-4-2J}/K)\cdot K\alpha\log (1/K\alpha)\\
&\ioe 2(N+1)/e.
\end{align*}

Enfin,
\begin{align*}
\sum_{j=0}^{J-1}\sum_{k=2}^{\alpha^{-1}}2\alpha \lfloor \alpha^{-4-2j}/k\rfloor &\ioe 2\alpha^{-3}\sum_{j=0}^{J-1}\alpha^{-2j}\cdot\sum_{k=2}^{\alpha^{-1}}1/k\\
&\ioe 2\alpha^{-3}\cdot \frac{\alpha^{-2J}-1}{\alpha^{-2}-1}\cdot \log(1/\alpha)\\
&\ioe \frac 83 \cdot \alpha\cdot (\alpha^{-4-2J} \cdot \alpha) \cdot \alpha\log(1/\alpha)\\
&\ioe \frac{8}{3e}\cdot \alpha\cdot (\alpha^{-4-2J}/K)\\
&\ioe \frac{8}{3e}\alpha(N+1).
\end{align*}

En réunissant ces trois majorations, on obtient bien que $G_{\alpha}(N)\ioe (N+1)(2/e+3\alpha)$.\fin

\smallskip

En particulier, si $\alpha^{-1} \soe 12$, on déduit de la proposition \ref{t43} que $\limsup_{x\vers \infty}G_{\alpha}(x)/x\ioe 1$.

\begin{prop}\label{t47}
Soit $J\in \Nat$ et $N=\alpha^{-4-2J}$. On a
$$
0 <\sum_{\alpha N\ioe n\ioe N}f_{\alpha}(n)/n < 1/N\alpha.
$$
\end{prop}
\dem

En effet,
\begin{align*}
\sum_{\alpha N\ioe n\ioe N}f_{\alpha}(n)/n&=\sum_{k=2}^{\alpha^{-1}}\Big (\frac{\alpha\lfloor N/k\rfloor}{\lfloor N/k\rfloor}-\frac{\alpha\lfloor N/k\rfloor}{\lfloor N/k\rfloor+1}\Big)\\
&=\alpha\sum_{k=2}^{\alpha^{-1}}\frac{1}{\lfloor N/k\rfloor+1}\\
&< \alpha N^{-1}\sum_{k=2}^{\alpha^{-1}}k\\
&\ioe 1/\alpha N.\fine
\end{align*}

\smallskip

Posons enfin
$$
H_{\alpha}(x)=\sum_{n\ioe x}f_{\alpha}(n)\{x/n\}.
$$
\begin{prop}\label{t48}
Soit $\alpha >0$ tel que $\alpha^{-1}$ soit un entier supérieur ou égal à $2$. On a
$$
\limsup_{x\vers \infty}|H_{\alpha}(x)/x| \soe \alpha \log (1/\alpha)+O(\alpha).
$$
\end{prop}
\dem

Soit $J\in \Nat$ et $N=\alpha^{-4-2J}$. On a
$$
H_{\alpha}(N)=\sum_{n\ioe N}f_{\alpha}(n)(N/n-\lfloor N/n\rfloor).
$$

D'après le principe de l'hyperbole, on a
\begin{align*}
\sum_{n\ioe N}f_{\alpha}(n)\lfloor N/n\rfloor &=\sum_{n< \alpha N}f_{\alpha}(n)\lfloor N/n\rfloor +\sum_{k\ioe \alpha^{-1}}F_{\alpha}(N/k)-F_{\alpha}(\alpha N-0)\alpha^{-1}\\
&=\sum_{n< \alpha N}f_{\alpha}(n)\lfloor N/n\rfloor +\alpha \sum_{2\ioe k\ioe \alpha^{-1}}\lfloor N/k\rfloor\expli{d'après \eqref{t45}}\\
&=\sum_{n< \alpha N}f_{\alpha}(n)\lfloor N/n\rfloor +\alpha N \big (\log (1/\alpha)+O(1)\big) +O(1).
\end{align*}

Par conséquent,
\begin{align*}
H_{\alpha}(N) &=\sum_{n< \alpha N}f_{\alpha}(n)\{ N/n\}+\sum_{\alpha N\ioe n\ioe N}f_{\alpha}(n)N/n-\alpha N \big (\log (1/\alpha)+O(1)\big) +O(1)\\
&=-\alpha N \big (\log (1/\alpha)+O(1)\big) +O(1/\alpha) \expli{d'après les propositions \ref{t43} et \ref{t47}}.
\end{align*}

On obtient le résultat annoncé en faisant tendre $J$ vers l'infini.\fin

\medskip

Finalement, la construction présentée montre que $c(\alpha) \soe \alpha \log (1/\alpha)+O(\alpha)$ si $\alpha^{-1}$ est un nombre entier supérieur ou égal à $12$. La réponse à la question \ref{t49} est bien négative.

\section{Questions ouvertes}

La principale question issue de ces remarques est celle de la détermination des constantes optimales $C_k$ (et dans une moindre mesure $D_k$) de la proposition \ref{t18}. Comme première question concrète, je propose la suivante. 
\begin{ques}
A-t-on
$$
|m_1(x)|\ioe x^{-3}\int_1^xt|M(t)|\,dt +O(1/x) \quad (x\soe 1)?
$$  
\end{ques}
Autrement dit, peut-on remplacer la constante $1,1$ du tableau final du \S\ref{t44} par 1?

\medskip

J'indique aussi deux questions qui apparaissent naturellement à propos du lemme d'Axer (cf. \S\S\ref{t38} et \ref{t41} ci-dessus). D'abord :

\begin{ques}
Quel est le comportement asymptotique de $c(\alpha)$ quand $\alpha$ tend vers $0$?
\end{ques}

\smallskip

Enfin soit $c_1(\alpha)$ la quantité définie comme $c(\alpha)$, mais en remplaçant la condition  
$$
\limsup_{x\vers \infty} G(x)/x \ioe 1
$$ 
par 
$$
\limsup_{n\vers \infty}|f(n)|\ioe 1.
$$

\begin{ques}
A-t-on 
$$
c_1(\alpha)=o\big(\alpha\log(1/\alpha)\big) \quad (\alpha \vers 0)?
$$
\end{ques}

\bigskip

\begin{center}

{\sc Remerciements}
\end{center}

 Je remercie le département de mathématiques de l'University College de Londres, et particulièrement Y. Petridis, pour d'excellentes conditions de travail lors de la préparation de cet article en juillet 2011. Je remercie également l'arbitre sollicité par l'Institut Steklov (et qui a tenu à rester anonyme), pour l'exemple présenté au \S\ref{t41} concernant le lemme d'Axer.


\providecommand{\bysame}{\leavevmode ---\ }
\providecommand{\og}{``}
\providecommand{\fg}{''}
\providecommand{\smfandname}{\&}
\providecommand{\smfedsname}{\'eds.}
\providecommand{\smfedname}{\'ed.}
\providecommand{\smfmastersthesisname}{M\'emoire}
\providecommand{\smfphdthesisname}{Th\`ese}

\bigskip

\footnotesize

\noindent Michel BALAZARD\\
Institut de Mathématiques de Luminy, FRE 3529\\
CNRS, Université d'Aix-Marseille\\
Campus de Luminy, Case 907\\
13288 Marseille Cedex 9\\
FRANCE\\
Adresse \'electronique : \texttt{balazard@iml.univ-mrs.fr}

\end{document}